\documentclass[12]{article}
\usepackage{epsfig}
\usepackage{amssymb,amsmath,amsfonts,soul,amsthm,enumerate}
\usepackage{ragged2e}
\usepackage{fancyhdr}
\usepackage{multirow}
\usepackage{hhline}
\usepackage{caption}
\usepackage{rotating}
\usepackage{graphicx}
\usepackage{mathrsfs}
\newtheorem{definition}{Definition}

\newtheorem{lemma}{Lemma}

\begin{document}
\textbf{MSC 11A41, 11Y40, 90C11}
\vskip 0.7cm
 \centerline {PROBABILISTIC METHOD OF PROVING TWIN PRIMES' INFINITUDE}
\vskip 0.3cm
\centerline{Nurlan N. Tashatov, Alua S. Turginbayeva, Serik A. Altynbek}
\vskip 0.2cm

\begin{abstract}
Statistical distribution of the primes in an arithmetic progression is considered. The estimation of prime numbers is given and combinatorial methods are used to calculate the twin primes on the available interval. The distribution and estimation of the number of primes on the twin primes rows are obtained. A new method of twin prime infinity is proposed.

\end{abstract}

\section{Introduction}
Zhang was submitted a paper to the Annals of Mathematics in 2013 which established a proof that states there are infinitely many pairs of prime numbers that differ by 70 million or less \cite{ZY:2014:AM}. Zhang's result set off a flurry of activity in the field.

It should be noted that a large number of original scientific papers, including \cite{RM:2013:Springer}-\cite{HR:2017:SMA}, are devoted to problems of prime numbers and of prime twin numbers. An analysis of these works shows that these works, as well as other works, give a tangible impetus and significant prerequisites for the further development of the theory of prime numbers.

Baibekov \cite{Baibekov:2016:SRR},\cite{Baibekov:2017:WAP},\cite{Baibekov:2009:FFD} propose a prime numbers matrix notion. In the present paper we proposed a new estimating method of the twin prime numbers in the prime numbers matrix, as well as a proof of twin prime infinity.

Consider the well-known asymptotic equality (known as the law of prime numbers distribution)
\[
\pi(x) \sim \frac{x}{\ln x}
\]

In 1850 Chebyshev proved that the ratio is bounded above and below by two positive constants near to $1$ for all $x$. This result became known as Chebyshev’s theorem
\[
0,89 \frac{x}{\ln x} \leq \pi(x)\leq 1,11\frac{x}{\ln x}
\]

Sylvester improved this estimation to
\[
  0,95 \frac{x}{\ln x}\leq \pi(x) \leq 1,05 \frac{x}{\ln x}
\]

\section{Estimate the number of primes in the interval}
Prime numbers theorem \cite{Nikitin:2010:ANT}, \cite{Balazar:2013:ALD}: we have for all $x\geq2$:
\begin{equation}\label {eq1}
\ln2 \cdot\frac{x}{\ln x}-2\leq \pi(x) \leq 4\cdot \ln2\cdot \frac{x}{\ln x}+\log_2x
\end{equation}

Taking into account the prime numbers matrix notation \cite{Baibekov:2016:SRR},\cite{Baibekov:2017:WAP},\cite{Baibekov:2009:FFD}, we consider an arithmetic progression of dimension in Baibekov's matrix $B_{p_3}=6$  (Table 1). It’s obvious that all prime numbers in the matrix are located on line 4 and line 6. Lines 1, 3 and 5 are always even, and numbers in line 2 is always divided by 3.
\begin{definition}
In the matrix, prime numbers pairs which satisfy condition $(5+6(k-1);1+6k)$ for any $k=\overline{(1..\infty)}$ will be called as twin numbers.
\end{definition}
The rows in which there are simple and composite numbers will be called simple strings. Lines 1, 2, 3, 5 are compound, and lines 4 and 6 are simple. We estimate the number of prime numbers in the interval of our matrix in rows with prime numbers. For this we take 4 and 6 lines.
\begin{table}[h!]
 \caption{Baibekov's matrix columns $B_{p_3}$.}
\begin{center}
\begin{tabular}{|c@{~\setlength{\arrayrulewidth}{1pt}\vline~}c|c|c@{~\setlength{\arrayrulewidth}{1pt}\vline~}c|c|c@{~\setlength{\arrayrulewidth}{1pt}\vline~}c|c|c@{~\setlength{\arrayrulewidth}{1pt}\vline~}c|c|c|c|c|c|}
            \hline
            &\begin{sideways}1 column\end{sideways}&\begin{sideways}2 column\end{sideways}&\begin{sideways}3 column\end{sideways}&\begin{sideways}4 column\end{sideways}&...&\begin{sideways}8 column\end{sideways}&\begin{sideways}9 column\end{sideways}&...&\begin{sideways}403 column\end{sideways}&\begin{sideways}404 column\end{sideways}&... &a &...&b&... \\ \cline {1-16}
             &\multicolumn{3}{c@{~\setlength{\arrayrulewidth}{1pt}\vline~}}{1interval}&\multicolumn{3}{c@{~\setlength{\arrayrulewidth}{1pt}\vline~}}{2 interval}&\multicolumn{3}{c@{~\setlength{\arrayrulewidth}{1pt}\vline~}}{3 interval}&\multicolumn{6}{c|}{4 interval}  \\ \cline {1-16}
            1& & & & & & & & & & & & & & &   \\ \cline {1-16}
             &2&8&14 &20&...&44&50&...&2414&2420&...&6a-4&...&6b-4&...  \\ \cline {1-16}
             &3&9&15 &21&...&45&51&...&2415&2421&...&6a-3&...&6b-3&...  \\ \cline {1-16}
             &4&10&16&22&...&46&52&...&2416&2422&...&6a-2&...&6b-2&...  \\ \cline {1-16}
             &5&11&17&23&...&47&53&...&2417&2423&...&6a-1&...&6b-1&... \\ \cline {1-16}
             &6&12&18&24&...&48&54&...&2418&2424&...&6a  &...&6b  &... \\ \cline {1-16}
             &7&13&19&25&...&49&55&...&2419&2425&...&6a+1&...&6b+1&... \\ \cline {1-16}
        \end{tabular}
    \end{center}
    \end{table}

Step 1. Using formula (\ref{eq1}), we take a lower estimate for all prime numbers to  $6b+1$
\[
\ln2\cdot\frac{6b+1}{\ln(6b+1)}-2
\]

Next we take, by the upper estimate of (\ref{eq1}), all prime numbers to  $6a-1$
\[
4\cdot \ln2\frac{6a-1}{\ln(6a-1)}+\log_2(6a-1)
\]

Thus, between the columns $(a,b)$, where $a\geq9$ and $(7a\leq b)$  or in the interval of the natural series $(6a-1,6b+1)$  there are at least
\[
\ln2\cdot\frac{6b+1}{\ln(6b+1)}-2-4\cdot \ln2\cdot\frac{6a-1}{\ln(6a-1)}-\log_2(6a-1)
\]
prime numbers.

Therefore, in lines 4 and 6 between the columns $(a, b)$, half of the estimated numbers are allocated
\begin{equation}\label{eq2}
\begin{gathered}
\frac{\ln2}{2}\cdot\frac{6b+1}{\ln(6b+1)}-1-2\cdot \ln2\cdot\frac{6a-1}{\ln(6a-1)}-\frac{\log_2(6a-1)}{2}
\end{gathered}
\end{equation}

Step 2. We split the columns in Table 1 with rounding into the whole part of the number by intervals according to the formula $ST_0=1$  and
\begin{equation}\label{eq3}
ST_n=e^{n!}
\end{equation}
where $n=\overline{1..\infty}$.
\begin{figure}[h]
\center\includegraphics[width=8cm]{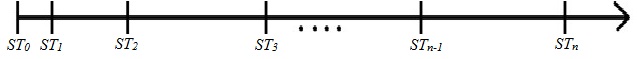}
\caption{The interval of partitioning the columns in the matrix.}
\end{figure}

Using formula (\ref{eq3}) for the table 1, we have the first interval $[ST_0;ST_1]$ will be between 1 and 3 columns, the second interval $(ST_1;ST_2]$ will be between 4 and 8 columns, the third $(ST_2;ST_3]$ interval will be between 9 and 403 columns, and the fourth $(ST_3;ST_4]$ interval will be between 404 and 26489122129 columns. And in each interval the number of columns increases exponentially according to the formula (\ref{eq3}).

The next task is to find the primes in each interval for each simple row.

For the following simplification of the calculations, we prove the following lemma
\begin{lemma}
For any number $i>2$, the following estimate of the number of primes in lines 4 and 6 in the interval of formula (\ref{eq2})
\end{lemma}
\begin{equation}\label{eq4}
\begin{gathered}
  \frac{\ln2}{2}\cdot\frac{6\cdot ST_i+1}{\ln(6\cdot ST_i+1)}-1-2\cdot \ln2\cdot\frac{6\cdot ST_{i-1}-1}{\ln(6\cdot ST_{i-1}-1)}-{} \\
  {}-\frac{\log_2(6\cdot ST_{i-1}-1)}{2}>\frac{ST_i}{\ln(ST_i)}-\frac{ST_{i-1}}{\ln(ST_{i-1})}.
\end{gathered}
\end{equation}
{\bf Proof}  For $i = 3, 4, 5$: formula (\ref{eq4}) takes the form
\[
3^{rd} \,\,\text{interval} - [int(e^{2!})+1; int(e^{3!})]:
\]
\[
88,046>63,543;
\]
\[
4^{th}\,\, \text{interval} - [int(e^{3!})+1; int(e^{4!})]:
\]
\[
2135665453,57>1103713354,83;
\]
\[
5^{nd}\,\, \text{interval} - [int(e^{4!})+1; int(e^{5!})]:
\]
\[
2,22672E50>1,08681E50.
\]

That the inequality.
Using (\ref{eq3}) and (\ref{eq4}) of the formula, we obtain the inequality:
\begin{equation}\label{eq5}
\begin{gathered}
\frac{\ln2}{2}\cdot\frac{6\cdot e^{i!}+1}{\ln(6\cdot e^{i!}+1)}-1-2\cdot \ln2\cdot\frac{6\cdot e^{(i-1)!}-1}{\ln(6\cdot e^{(i-1)!}-1)}-\frac{\log_2(6\cdot e^{(i-1)!}-1)}{2}> {} \\
{} >\frac{e^{i!}}{\ln(e^{i!})}-\frac{e^{(i-1)!}}{\ln(e^{(i-1)!})};
\end{gathered}
\end{equation}

\begin{equation*}
\begin{gathered}
\frac{\ln2}{2}\cdot\frac{6\cdot e^{i!}+1}{\ln(6\cdot e^{i!}+1)}-\frac{e^{i!}}{\ln(e^{i!})}>1+2\cdot \ln2\cdot
\frac{6\cdot e^{(i-1)!}-1}{\ln(6\cdot e^{(i-1)!}-1)}+{}\\{}+\frac{\log_2(6\cdot e^{(i-1)!}-1)}{2}-\frac{e^{(i-1)!}}{\ln(e^{(i-1)!})}.
\end{gathered}
\end{equation*}

In the formula (\ref{eq5}) we apply arithmetic operations, and we strengthen the inequality. The sequence of calculations is shown below:
for $i>5$
\[
\frac{\ln2}{2}\cdot \frac{6\cdot e^{i!}}{\ln(6\cdot e^{i!}+1)}\geq\frac{\ln2}{2}\cdot \frac{6\cdot e^{i!}}{\ln(6\cdot e^{i!})}-1;\,\,\,\, \frac{e^{i!}}{\ln(e^{i!})}<1,009\cdot \frac{e^{i!}}{\ln(6\cdot e^{i!})};
\]
\[
\ln2\cdot\frac{e^{(i-1)!}}{(i-1)!}>\frac{(i-1)!}{2\cdot \ln2}; \,\,\,\, \frac{e^{(i-1)!}}{{(i-1)!}^i}\cdot i>1;
\]

\begin{eqnarray*}
\begin{gathered}
\frac{\ln2}{2}\cdot \frac{6\cdot e^{i!}}{\ln(6\cdot e^{i!})}-\frac{e^{i!}}{\ln(e^{i!})}-1>1+2\cdot \ln2\cdot\frac{6\cdot e^{(i-1)!}}{\ln(6\cdot e^{(e-1)!})}+\frac{\ln(6\cdot e^{(i-1)!})}{2\cdot \ln2};
\end{gathered}
\end{eqnarray*}
\begin{eqnarray*}
\begin{gathered}
\bigg(\frac{\ln2\cdot 6}{2}-1,009\bigg)\cdot\frac{e^{i!}}{\ln(6\cdot e^{i!})}>2+\ln2\cdot \frac{6\cdot e^{(i-1)!}}{\ln(6\cdot e^{(i-1)!})}+\frac{\ln6}{2\cdot \ln2}+\frac{\ln(e^{(i-1)!})}{2\cdot \ln2};
\end{gathered}
\end{eqnarray*}
\begin{eqnarray*}
\begin{gathered}
2\cdot\frac{e^{i!}}{\ln(6\cdot e^{i!})}>2+\ln2\cdot \frac{6\cdot e^{(i-1)!}}{(i-1)!}+\frac{\ln6}{2\cdot \ln2}+\frac{(i-1)!}{2\cdot \ln2};
\end{gathered}
\end{eqnarray*}
\begin{equation}\label{eq6}
  2\cdot\frac{e^{i!}}{\ln(6\cdot e^{i!})}>3+\ln2\cdot\frac{7\cdot e^{(i-1)!}}{(i-1)!}.
\end{equation}

Further, inequality (\ref{eq6}) for $i>5$  by induction starts to prove.

Assume for $i=t$, inequality (\ref{eq6}) is true, that is, it has the form
\[
2\cdot\frac{e^{t!}}{\ln(6\cdot e^{t!})}>3+\ln2\cdot\frac{7\cdot e^{(t-1)!}}{(t-1)!}
\]

It is necessary to prove the inequality (\ref{eq6}) for $i=t+1$. Substituting $i=t+1$ in formula (\ref{eq6}) and applying arithmetic operations, we obtain the following inequality
\begin{equation*}
\begin{gathered}
  2\cdot\frac{(e^{t!})^{t+1}}{\ln(6\cdot e^{(t+1)!})}>2\bigg[\frac{e^{t!}}{\ln(6\cdot e^{t!})}\bigg]^{t+1}>3+7\cdot \ln2\cdot\frac{(e^{(t-1)!})^t}{(t-1)!t};
\end{gathered}
\end{equation*}
\[
\frac{e^{t!}}{\ln(6\cdot e^{t!})}>\frac{3}{2}+\frac{7}{2}\cdot\ln2\cdot\frac{ e^{(t-1)!}}{(t-1)!};
\]
\begin{equation*}
\begin{gathered}
  \bigg[\frac{3}{2}+\frac{7}{2}\cdot\ln2\cdot\frac{e^{(t-1)!}}{(t-1)!}\bigg]^{t+1}> {\bigg(\frac{3}{2}\bigg)}^{t+1}+\bigg(\frac{7}{2}\ln2\bigg)^{t+1}\cdot\frac{e^{(t-1)!\cdot t+(t-1)!}}{{(t-1)!}^{t+1}};
\end{gathered}
\end{equation*}
\begin{equation*}
\begin{gathered}
  \bigg[\frac{3}{2}+\frac{7}{2}\cdot\ln2\cdot\frac{e^{(t-1)!}}{(t-1)!}\bigg]^{t+1}>
  {\bigg(\frac{3}{2}\bigg)}^{t+1}+\bigg(\frac{7}{2}\ln2\bigg)^{t+1}\cdot\frac{e^{(t-1)!\cdot t}}{{(t-1)!}\cdot t}\cdot\frac{e^{(t-1)!}}{{(t-1)!}^t}\cdot t;
\end{gathered}
\end{equation*}
\begin{equation*}
\begin{gathered}
     {\bigg(\frac{3}{2}\bigg)}^{t+1}+\bigg(\frac{7}{2}\ln2\bigg)^{t+1}\cdot\frac{e^{(t-1)!\cdot t}}{{(t-1)!}\cdot t}\cdot\frac{e^{(t-1)!}}{{(t-1)!}^t}\cdot t>
     3+7\cdot\ln2\cdot\frac{e^{t!}}{t!}.
\end{gathered}
\end{equation*}
Thus, according to the method of mathematical induction, inequality (\ref{eq6}) is valid for any natural number $i>2$.

\textbf{The lemma is proved.}

\section{Estimation of prime numbers of twins by the method of mathematical statistics and probability theory}
Using the elements of probability theory and combinatorics, the proof of the infinity of the numbers of twins.

Thus, using Lemma 1 for the convenience of calculations, we will use an underestimated estimate the number of primes in each interval for $i>2$ for 4 and 6 lines
\begin{equation}\label{eq7}
  \frac{e^{i!}}{\ln(e^{i!})}-\frac{e^{(i-1)!}}{\ln(e^{(i-1)!})}
\end{equation}

Taking any number from the interval (numbers are taken from 4 and 6 lines), and applying formulas (\ref{eq2}) and (\ref{eq4}), we can calculate the probability that this number is prime or not according to the classical probability formula \cite{Gmurman:2018:TBMS}.
\begin{equation}\label{eq8}
\begin{gathered}
  P \left(\frac{number\;of\;prime\;numbers}{number\;of\;
  columns\;in\;interval}\right)=\frac{\frac{e^{i!}}{\ln(e^{i!})}-\frac{e^{(i-1)!}}{\ln(e^{(i-1)!})}}{e^{i!}-e^{(i-1)!}}.
\end{gathered}
\end{equation}

Now calculate the probability of occurrence of the twin primes in each interval, using the definition (1). This problem follows from the probability multiplication theorem \cite{Gmurman:2018:TBMS} for the independent events $P(AB)=P(A)\cdot P(B)$

The solution of the problem: $m_1,m_2$  -- the number of primes in simple intervals of the interval is calculated from the formula (\ref{eq2}). We can objectively assume that the number of primes in the interval $m_1 \approx m_2$, since in the future we will use the underestimated density of estimating primes by formula (\ref{eq4}) [13].
\begin{equation}\label{eq9}
P\left(\frac{m_1}{n}\right)\cdot P\left(\frac{m_2}{n}\right)=\left(\frac{\frac{e^{i!}}{\ln(e^{i!})}-\frac{e^{(i-1)!}}{\ln(e^{(i-1)!})}}{e^{i!}-e^{(i-1)!}}\right)^2
\end{equation}

Using formula (\ref{eq9}) and the number of columns in the interval, we can estimate the mathematical expectation of the discrete distribution of the numbers of twins in each interval.
\begin{equation}\label{eq10}
\begin{gathered}
  M_i(X)=\left(\frac{\frac{e^{i!}}{\ln(e^{i!})}-\frac{e^{(i-1)!}}{\ln(e^{(i-1)!})}}{e^{i!}-e^{(i-1)!}}\right)^2\cdot (e^{i!}-e^{(i-1)!})  =\frac{\left(\frac{e^{i!}}{\ln(e^{i!})}-\frac{e^{(i-1)!}}{\ln(e^{(i-1)!})}\right)^2}{e^{i!}-e^{(i-1)!}}
\end{gathered}
\end{equation}
where $\left(\frac{\frac{e^{i!}}{\ln(e^{i!})}-\frac{e^{(i-1)!}}{\ln(e^{(i-1)!})}}{e^{i!}-e^{(i-1)!}}\right)^2$ probability of the appearance of a number of twins satisfying the definition (1) in the interval $(e^{(i-1)!}, e^{i!}]$, and $(e^{i!}-e^{(i-1)!})$ is the number of columns between the interval.

Thus, for each interval in Fig.1, we estimated the mathematical expectation of the appearance of the twin primes.

The mathematical expectation of the number of identical positions.

The required number will be calculated for the case of coincidence of the positions of the supposed pairs of primes in both rows. The case where the positions are relatively shifted by one is easily reduced to the one considered. We will also consider only the positions on which, in principle, there can be a prime number.

Let there be two lines of length $n$, for which $m_1$ and $m_2$ coordinates, respectively, with values of one, the rest - with values of zero. Suppose that $(m_1\sim m_2)$ and $m_1+m_2<<n$, that is, units in the sum is quite small. We calculate the number of cases in which the same exactly $s$ positions of units in both rows. Let $m_1$ units of the first vector somehow located, and such different arrangements can be $\binom{n}{m_1}$. From these positions of $s$ positions, one can choose  $\binom{m_1}{s}$  in ways that must be combined with any of  $\binom{n-m_1}{m_2-s}$ locations of the remaining units of the second row at the positions of zeros of the first. Thus, the number of locations of units in both rows, provided that exactly $s$ positions are the same, is equal to $T(s)=\binom{n}{m_1} \binom{m_1}{s} \binom{n-m_1}{m_2-s}$. Then
\[
\sum_{s=0}^{s=m_1}T(s)=\binom{n}{m_1}\sum_{s=0}^{s=m_1}\binom{m_1}{s} \binom{n-m_1}{m_2-s}=\binom{n}{m_1}\binom{n}{m_2}
\]

The mathematical expectation of the number of identical positions under the condition of a random arrangement of units in both rows is equal, by definition,
\[
\frac{1}{\binom{n}{m_1}\binom{n}{m_2}}\sum_{s=0}^{s=m_1}sT(s)=\frac{1}{\binom{n}{m_2}}\sum_{s=0}^{s=m_1}s\binom{m_1}{s} \binom{n-m_1}{m_2-s}=\frac{m_1}{\binom{n}{m_2}}\sum_{s=0}^{s=m_1}\binom{m_1-1}{s-1} \binom{n-m_1}{m_2-s}=\frac{m_1 m_2}{n}
\]
where
\[
n=e^{i!}-e^{(i-1)!}, \,\,\,m_1\approx m_2=\frac{e^{i!}}{\ln(e^{i!})}-\frac{e^{(i-!)!}}{\ln(e^{(i-!)!})}
\]

In formula (\ref{eq10}), summing up all the mathematical expectations in each interval, we estimate the series
\begin{equation}\label{eq11}
  \sum_{i=3}^{n}M_i(X)=\sum_{i=3}^{n}\frac{\left(\frac{e^{i!}}{\ln(e^{i!})}-\frac{e^{(i-1)!}}{\ln(e^{(i-1)!})}\right)^2}{e^{i!}-e^{(i-1)!}}
\end{equation}

To analyze the convergence of the series, we apply the verification of d'Alembert to formula (\ref{eq11})
\begin{equation}\label{eq12}
  \frac{a_n}{a_{n-1}}\to ?
\end{equation}

The calculations of formula (\ref{eq12}) are given below
\begin{equation}\label{eq13}
\begin{gathered}
  \frac{a_n}{a_{n-1}}=\frac{\Big(\frac{e^{n!}}{\ln(e^{n!})}-\frac{e^{(n-1)!}}{\ln(e^{(n-1)!})}\Big)^2}{e^{n!}-e^{(n-1)!}}
  \vdots  \frac{\Big(\frac{e^{(n-1)!}}{\ln(e^{(n-1)!})}-\frac{e^{(n-2)!}}{\ln(e^{(n-2)!})}\Big)^2}{e^{(n-1)!}-e^{(n-2)!}}={} \\
   {}
   =\frac{[e^{n!}\cdot(n-1)!-e^{(n-1)!}\cdot n!]^2}{(n!\cdot(n-1)!)^2\cdot(e^{n!}-e^{(n-1)!})}\cdot \frac{((n-1)!\cdot(n-2)!)^2\cdot(e^{(n-1)!}-e^{(n-2)!})}{[e^{(n-1)!}\cdot(n-2)!-e^{(n-2)!}\cdot(n-1)!]^2}={}
   \\
   {}
   =\frac{[e^{n!}-e^{(n-1)!}\cdot n]^2}{n^2\cdot(e^{n!}-e^{(n-1)!})}\cdot \frac{(e^{(n-1)!}-e^{(n-2)!})}{[e^{(n-1)!}-e^{(n-2)!}\cdot(n-1)]^2}
\end{gathered}
\end{equation}

In formula (\ref{eq13}) expanding the brackets and canceling the numerator and denominator of $e^{(2\cdot n!+(n-1)!)}$, and putting $n\to\infty$, we obtain $\frac{a_n}{a_{(n-1)}} =\infty$. This implies that the series (\ref{eq11}) diverges. Thus, the mathematical expectation from formula (\ref{eq11}) shows that the number of twin primes tends to infinity.

\section{Conclusion}
In this paper, the Baibekov's matrix notation is used for investigation of the prime twin numbers distribution. Then in the paper we prove a number of lemmas and theorems with the help of which and the Dirichlet and Euler theorems are proposed to prove the infinity of the number of prime twins.

It is obtained that the density of prime numbers in the process of transition to the region of large natural numbers continuously in simple lines falls. In the first intervals, prime numbers often occur, which are more closely aligned. As they move into the region of the following intervals, they are less common. But in contrast to the number of columns in subsequent intervals in the exponential order increases. Thus, the infinity of prime twin numbers is proved. We carried out a numerical experiment on program C\#, calculated the number of prime numbers for each particular interval.\newline
\textit{{\bf{Acknowledgement.}} The present work was carried out on the basis of the scientific and technical program {N\textsuperscript{\underline{o}}} BR 05236075 of the Science Committee of the Ministry of Education and Science of the Republic of Kazakhstan on the priority development of science "Mangilik El"}

{
Associate professor Nurlan Tashatov, doctoral student Alua Turginbayeva,
L.N. Gumilyov Eurasian National University,
st. Pushkin 11,
010000 Astana, Kazakhstan,
emails: tash.nur@mail.ru, ssalua@mail.ru,
PhD Serik Altynbek,
RSE "Institute of Information and Computing Technologies" of the Science Committee of the Ministry of Education and Science of the Republic Kazakhstan,
st. Kabanbay batyr 8,
010000 Astana, Kazakhstan,
email: serik\_aa@bk.ru}



\begin{thebibliography}{999}
\let\titem\relax 

\bibitem{ZY:2014:AM}
Zhang Y.,
\newblock Bounded Gaps between Primes.
\newblock{\titem Annals of Mathematics.} 2014. pp. 1121--1174. DOI:https://doi.org/10.4007/annals.2014.179.3.7

\bibitem{RM:2013:Springer}
Melnik R., Kotsireas I.,
\newblock Advances in Applied Mathematics, Modeling, and Computational Science,
\newblock{\titem Springer Science+Business Media, New York, U.S.A.}, (2013). MR2963948.

\bibitem{HM:2015:MMGH}
Muller H.,
\newblock A Journey Through the World of Primes,
\newblock {\titem Mitt. Math. Ges. Hamb.}, 35(2015), pp. 5–-18. Zbl 1343.11007.

\bibitem{MO:2014:AMS}
Overholt M.,
\newblock A Course in Analytic Number Theory,
\newblock {\titem American Mathematical Society}, 2014. MR3290245. Zbl 1335.11003.

\bibitem{Ribenboim:1996:Springer}
Ribenboim P.,
\newblock The New Book of Prime Number Records,
\newblock {\titem Springer-Verlag, New York, U.S.A., 3rd  ed.}, 1996. MR1377060. Zbl 0856.11001.

\bibitem{HR:2017:SMA}
Rezgui H.,
\newblock Conjecture of twin primes (still  unsolved problem in number theory) an expository essay,
\newblock{\titem Surveys in Mathematics and its Applications}, Volume 12 (2017), pp. 229--252.

\bibitem{Baibekov:2016:SRR}
 Baibekov S.N., Durmagambetov A.A.,
  \newblock Infinite Number of Twin Primes,
 \newblock {\titem Advances in Pure Mathematics}.  2016, no.~6, pp. 954--971.

\bibitem{Baibekov:2017:WAP}
Baibekov S.N., Durmagambetov A.A.,
\newblock Proof of an Infinite Number of Primes-Twins,
\newblock {\titem International Journal of Scientific and Innovative Mathematical Research (IJSIMR)}, 2017, 5 (19), pp. 6--17.

\bibitem{Baibekov:2009:FFD}
Baibekov S.N., Dossayeva A.A.,
\newblock Development of the Matrix of Primes and Proof of an Infinite Number of Primes-Twins.,
\newblock https://arxiv.org/ftp/arxiv/papers/1805/1805.00346.pdf.

\bibitem{Nikitin:2010:ANT}
Nikitin N.D.,
\newblock Algebra i teoriya chisel: Uchebnoe posobie,
\newblock {\titem Algebra and Number Theory: Tutorial} Penza, PGPU, 2010. p. 96.

\bibitem{Balazar:2013:ALD}
Balazar M.,
\newblock Asimptoticheskiy zakon raspredeleniya prostyh chisel,
\newblock {\titem Asymptotic Law of Distribution of Prime Numbers}. Moscow, MTSNMO, 2013. p. 64.

\bibitem{Gmurman:2018:TBMS}
Gmurman V.E.,
\newblock Teoriya veroyatnostey I matematicheskaya statistika. Uchebnik,
\newblock {\titem Theory of Probability and Mathematical Statistics. Tutorial}, Moscow, 2018. p. 480.

\bibitem{Altynbek:2018:BPU}
Abdymanapov S.A., Altynbek S.A., Turginbayeva A.S.,
\newblock Calculating the number of twin primes with specified distance between them based on the simplest probabilistic model,
\newblock {\titem University Politehnica of Bucharest Scientific Bulletin-Series A-Applied Mathematics and Physics}, Vol. 80, Iss.3, 2018. pp. 183-190.

\end{thebibliography}
\end{document}